\documentclass[12pt,a4paper,reqno]{amsart}
\usepackage{amssymb,bm}
\usepackage{latexsym}
\usepackage{exscale}
\usepackage{mathrsfs}
\usepackage{color}

\newcommand{\sG}{\mathscr G}
\newcommand{\sB}{\mathscr B}

\newcommand{\tri}[1]{|\!|\!| #1 |\!|\!|}

\newcommand{\fpq}{{\mathfrak{f}^q_{1}}}

\newcommand{\ba}{{\mathbf a}}
\newcommand{\be}{{\mathbf e}}
\newcommand{\ben}{{\be^*_n}}

\newcommand{\bL}{{\mathbf L}_N}
\newcommand{\tL}{{\mathbf{\widetilde{L}}}_N}
\newcommand{\bq}{{\mathbf q}}
\newcommand{\bs}{{\mathbf s}}

\newcommand{\bfe}{{\boldsymbol\e}}
\newcommand{\bfeta}{{\boldsymbol\eta}}

\newcommand{\cD}{{\mathcal{D}}}

\newcommand {\SC} {{\mathbb C}}

\newcommand {\SK} {{\mathbb K}}
\newcommand {\SN} {{\mathbb N}}
\newcommand {\SR} {{\mathbb R}}
\newcommand {\ST} {{\mathbb T}}

\newcommand {\SX} {{\mathbb X}}

\newcommand {\SZ} {{\mathbb Z}}

\newcommand {\tmu} {{\tilde\mu}}
\newcommand {\tg} {{\tilde{g}}}

\newcommand {\al} {{\alpha}}

\newcommand {\e} {{\varepsilon}}
\newcommand {\ga} {{\gamma}}
\newcommand {\Ga} {{\Gamma}}
\newcommand {\tGa} {{\widetilde\Gamma}}
\newcommand {\tsigma} {{\widetilde\sigma}}
\newcommand {\la} {{\lambda}}
\newcommand {\La} {{\Lambda}}

\newcommand {\chii} {{\chi^{(1)}_I}}

\newcommand {\hg} {{\hat{g}}}
\newcommand {\vg} {{\check{g}}}

\newcommand{\ttK}{\texttt{K}}

\def\supp{\mathop{\rm supp}}
\def\dist{\mathop{\rm dist}}
\def\sig{\mathop{\rm sign}}
\def\co{\mathop{\rm co}}

\numberwithin{equation}{section}
\newtheorem{theorem}{Theorem}[section]
\newtheorem{lemma}[theorem]{Lemma}

\newtheorem{corollary}[theorem]{Corollary}
\newtheorem{Remark}[theorem]{Remark}
\newtheorem{proposition}[theorem]{Proposition}

\newtheorem{example}[theorem]{Example}

\newcommand {\Proofof}[1] {\noindent{\bf P{\footnotesize\bf ROOF} of {#1}: } \ }

\newcommand {\Proof} {\noindent{\bf P{\footnotesize\bf ROOF}: } \ }
\newcommand {\ProofEnd} {
             \begin{flushright} \vskip -0.2in $\Box$ \end{flushright}}

\newcommand{\Ba}[1]{\begin{array}{#1}}
\newcommand{\Ea}{\end{array}}
\newcommand{\Bd}{\begin{description}}
\newcommand{\Ed}{\end{description}}
\newcommand{\Be}{\begin{equation}}
\newcommand{\Ee}{\end{equation}}
\newcommand{\Bea}{\begin{eqnarray}}
\newcommand{\Eea}{\end{eqnarray}}
\newcommand{\Beas}{\begin{eqnarray*}}
\newcommand{\Eeas}{\end{eqnarray*}}
\newcommand{\Benu}{\begin{enumerate}}
\newcommand{\Eenu}{\end{enumerate}}
\newcommand{\Bi}{\begin{itemize}}
\newcommand{\Ei}{\end{itemize}}

\newcommand{\BR}{\begin{Remark} \em}
\newcommand{\ER}{\end{Remark}}
\newcommand{\BE}{\begin{example} \em}
\newcommand{\EE}{\end{example}}

\newcommand {\mand} {{\quad\mbox{and}\quad}}
\renewcommand {\mid} {{\,\,\,\colon\,\,\,}}

\newcommand{\bline}{{\bigskip

\noindent}}

\newcommand{\sline}{{\smallskip

\noindent}}

\newcommand {\bone} {{\bf 1}}

\newcommand {\tG} {{\widetilde \Ga}}

\newcommand{\cupdot}{\mathbin{\mathaccent\cdot\cup}}

\newcounter{reg}
\begin{document}

\title[Lebesgue inequalities in general bases]{Lebesgue inequalities for the greedy algorithm in general bases}%
\author{Pablo M. Bern\'a}
\address{Pablo M. Bern\'a
\\
Instituto Universitario de Matem\'atica Pura y Aplicada
\\
Universitat Polit\`ecnica de Val\`encia
\\
46022 Valencia, Spain} \email{pmbl1991@gmail.com}

\author{\'Oscar Blasco}
\address{\'Oscar Blasco
\\
Departamento de An\'alisis Matem\'atico
\\
Universidad de Valencia, Campus de Burjassot
\\
46100 Valencia, Spain} \email{oscar.blasco@uv.es}

\author{Gustavo Garrig\'os}
\address{Gustavo Garrig\'os
\\
Departamento de Matem\'aticas
\\
Universidad de Murcia
\\
30100 Murcia, Spain} \email{gustavo.garrigos@um.es}

\thanks{First author partially supported by GVA PROMETEOII/2013/013 and grant 19368/PI/14 below.
Second author partially supported by grant MTM2014-53009-P (MINECO, Spain).
Third author partially supported by grants MTM2013-40945-P, MTM2014-57838-C2-1-P (MINECO, Spain) and 19368/PI/14 (\emph{Fundaci\'on S\'eneca}, Regi\'on de Murcia, Spain).}



\begin{abstract}
We present various estimates for the Lebesgue constants of the thresholding greedy algorithm, in the case of general bases in Banach spaces.
We show the optimality of these estimates in some situations.
Our results recover and slightly improve various estimates appearing earlier in the literature.
\end{abstract}

\maketitle

\section{Introduction }
\setcounter{equation}{0}\setcounter{footnote}{0}
\setcounter{figure}{0}

Let $\SX$ be a Banach space (over $\SK=\SR$ or $\SC$) and $\{\be_n, \ben\}_{n=1}^\infty$ a biorthogonal system such that $\sB=\{\be_n\}$ has dense span in $\SX$ and $0<\kappa_1\leq \|\be_n\|,\|\ben\|\leq \kappa_2<\infty$. Examples include (semi-normalized) Schauder bases $\sB$, as well as more general structures (such as Markushevich bases \cite{Hajek}). As suggested in \cite{Wo,Wo3}, greedy algorithms can be considered in this generality, by formally associating with every $x\in\SX$ the series $x\sim\sum_{n=1}^\infty e^*_n(x)e_n$.
Note that $\lim_{n\to\infty}\ben(x)=0$, so one may speak of decreasing rearrangements of $\{\ben(x)\}$.

We recall a few standard notions about greedy algorithms;  see e.g.
\cite{Tem1,Tem15} for a detailed presentation and background.
We say that a finite set $\Ga\subset\SN$ is a greedy set for $x\in\SX$, denoted $\Ga\in\sG(x)$, if
\[
\min_{n\in \Ga}|\ben(x)|\,\geq\,\max_{n\in\Ga^c}|\ben(x)|,
\]
and write $\Ga\in\sG(x,N)$ if in addition $|\Ga|=N$. A \emph{greedy operator of order $N$} is a mapping  $G:\SX\to\SX$ such that
\[
Gx=\sum_{n\in \Ga_x}\ben(x)\be_n, \quad \mbox{for some }\Ga_x\in \sG(x,N).
\]
We write $\sG_N$ for the set of all greedy operators of order $N$, and $\sG=\cup_{N\geq1}\sG_N$. Given $G,G'\in \sG$ we shall write
$G'<G$ whenever $G\in \sG_N$ and $G'\in \sG_M$ with $M<N$ and $\Gamma'_x\subset \Gamma_x$.

Likewise, for every \emph{finite} set $A\subset \SN$ we consider the projection operator
\[P_Ax=\sum_{n\in A}\ben(x)\be_n,\]
and the ``complement'' projection $P_{A^c}=I-P_A$.

\

Greedy operators are frequently used for $N$-term approximation. As usual, we let $\Sigma_N=\big\{\sum_{A}a_n\be_n\mid |A|\leq N,\;a_n\in\SK\big\}
$ 
and
$\sigma_N(x)=\dist (x,\Sigma_N)$. 
To quantify the efficiency of greedy approximation one defines, for each $N=1,2,\ldots$, the smallest number $\bL$ such that\Be
\|x-Gx\|\leq \bL\,\sigma_N(x), \quad \forall\;x\in\SX,\;
\; \forall\;G\in\sG_N.
\label{bL}\Ee
This is sometimes called a Lebesgue-type inequality for the greedy algorithm \cite{Tem15}, and $\bL$ is its associated Lebesgue-type constant. Likewise, one may consider ``expansional'' $N$-term approximations and $\widetilde\sigma_N(x)=\inf\{\|x-P_Ax\|\mid |A|\leq N\}$, and define the smallest $\tL$ such that
\Be
\|x-Gx\|\leq \tL\,\widetilde\sigma_N(x), \quad \forall\;x\in\SX,\;
\; \forall\;G\in\sG_N.
\label{tL}\Ee

\

A celebrated result of Konyagin and Temlyakov \cite{KT} establishes that $\bL=O(1)$ if and only if $\sB$ is unconditional and democratic.
Explicit estimates for $\bL$ have been obtained in various contexts for greedy bases \cite{Wo3,AlWo,DKOSZ}, quasi-greedy bases \cite{TYY2,DST,GHO,DKO,AA1}, and a few examples of non quasi-greedy bases \cite{Tem98trig,Tem98haar,Osw}. The goal of this paper is to present these inequalities in a more general setting, and improve them as much as possible so that they actually become optimal
in certain Banach spaces. This of course depends on the quantities used for the bounds, which we list next.

\Bi
\item  Unconditionality constants: \[k_N=\sup_{|A|\leq N}\|P_A\|\mand k^c_N=\sup_{|A|\leq N}\,\|I-P_A\|.\]

\item  Quasi-greedy constants\footnote{We use the notation $\|G\|=\sup_{x\not=0}\|Gx\|/\|x\|$, even if $G:\SX\to\SX$ is a non-linear map.}:
\[g_N=\sup_{G\in\cup_{k\leq N}\sG_k}\,\|G\|\  \mand  g^c_N=\sup_{G\in\cup_{k\leq N}\sG_k}\,\|I-G\| .
\]
We shall also use
\[\hg_N=\min\{g_N,g^c_N\}\mand \tilde g_N= \sup_{G\in\cup_{k\leq N}\sG_k, G'<G}\,\|G-G'\|\
.\]

\

\item Democracy (and superdemocracy) constants:
\[
\mu_N=\sup_{|A|=|B|\leq N}\frac{\|\bone_A\|}{\|\bone_B\|}\mand
\tmu_N=\sup_{{|A|=|B|\leq N}\atop{\bfe,\bfeta\in\Upsilon}}\frac{\|\bone_{\bfe A}\|}{\|\bone_{\bfeta B}\|},
\]
and their counterparts for disjoint sets given by
\[
\mu^d_N=\sup_{{|A|=|B|\leq N}\atop{A\cap B=\emptyset}}\frac{\|\bone_A\|}{\|\bone_B\|}\mand
\tmu^d_N=\sup_{{{|A|=|B|\leq N}\atop{A\cap B=\emptyset}}\atop{\bfe,\bfeta\in\Upsilon}}\frac{\|\bone_{\bfe A}\|}{\|\bone_{\bfeta B}\|}
\]
\item A-property constants:
\[
\nu_N =\sup\left\{\frac{\|\bone_{\bfe A}+x\|}{\|\bone_{\bfeta B}+x\|}\mid \mbox{ {\footnotesize $|A|=|B|\leq N, \;\bfe,\bfeta\in\Upsilon,\;  |x|_\infty\leq 1,\; A\cupdot B\cupdot \mbox {\normalsize$x$}$} } \right\}.
\]
\Ei

\

We are using the standard notation
\[\bone_A=\sum_{n\in A}\be_n\mand \bone_{\bfe A}= \sum_{n\in A}\e_n\be_n,\quad {\rm if}\;\bfe=\{\e_n\}.\]
Here $\bfe=\{\e_n\}\in\Upsilon$ means that $|\e_n|=1$ for all $n$ (where $\e_n$ could be real or complex). We also set $|x|_\infty=\sup_n|\ben(x)|$ and $\supp x=\{n\mid \ben(x)\not=0\}$, and we write $A\cupdot B\cupdot x$ to mean that $A,B$ and $\supp x$ are pairwise disjoint.

All these are natural constants in the greedy literature, and often it is not hard to compute them explicitly; see  $\S\ref{examples}$ below for some examples. Let us point out some elementary inequalities for the less frequent constants $\tg_N$ and $\nu_N$.
\BR\label{Re1} For each $N\in \mathbb N$ we have \Be g_N\leq \tilde g_N\le \min\{2\hg_N,\, g_N g_N^c,\,k_N\}.\Ee
Indeed,  $g_N\leq \tilde g_N\le k_N$ is obvious by definition and $\tilde g_N\le 2\hg_N$ follows easily from the triangle inequality. Finally, for each $G\in \cup_{k\leq N}\sG_k$ and $G'<G$ we can write $Gx-G'x= \sum_{n\in \Gamma_x\setminus \Gamma'_x} e_n^*(x)e_n$ with $\Gamma_x\setminus \Gamma'_x\in \cup_{k\leq N}\sG(x-G'x, k)$; hence
$$\|Gx-G'x\|\le g_N \|x-G'x\|\le g_Ng_N^c \|x\|.$$
\ER
\BR \label{Re2}For each $N\in \mathbb N$ we have
\Be \max\{\tmu^d_N, \mu_N\}\le \nu_N\le g_N^c+ g_N\tmu^d_N.\Ee
Indeed, the inequalities $\tmu^d_N\le \nu_N$ and $\mu_N\le \nu_N$ follow selecting $x=0$ and $x=\bone_{A\cap B}$ respectively in the definition of $\nu_N$.
On the other hand, for each $|A|=|B|\leq N, \;\bfe,\bfeta\in\Upsilon,\;  |x|_\infty\leq 1,\; A\cupdot B\cupdot \mbox {\normalsize$x$}$ we have
$\|x\|\le g_N^c \|\bone_{\bfe B}+x\|$ and
$\|\bone_{\bfe A}\|\le \tmu^d_N\|\bone_{\bfe B}\|\le \tmu^d_Ng_N\|\bone_{\bfe B}+x\|$. Hence
the inequality $\nu_N\le g_N^c+ g_N\tmu^d_N$ is easily obtained.
\ER

The above mentioned constants are also natural \emph{lower} bounds for the Lebesgue inequalities.

\begin{proposition}\label{P1}
For all $N\geq1$ we have
\Be
\bL\,\geq\,\max\big\{ k^c_N\,,\;\tL\big\},\mand \tL\,\geq\,\max\big\{ g^c_N\,,\;\nu_N\,,\;\mu_N\;,\tfrac1{2\kappa}\tmu_N\;\big\},
\label{lower}\Ee
with $\kappa=1$ for real spaces, and $\kappa=2$ for complex spaces.
\end{proposition}

\

We shall present two results concerning upper bounds.

\begin{theorem}\label{th2}
For all $N\geq1$ we have
\Be\label{LN2}
\bL \leq k^c_{2N}\,\nu_N\mand \tL \leq g^c_{N}\,\nu_N.
\Ee
Moreover, there exists $(\SX,\sB)$ for which both equalities are attained.
\end{theorem}

\begin{theorem}\label{th3}
For all $N\geq1$ we have
\Be
\bL \leq k^c_{2N}\,+\,\tilde g_N\,\tmu_N\;\mand\;
\tL \leq g^c_{N}\,+\,\tilde g_N\,\tmu_N.
\label{LN3}\Ee
Moreover, there exists $(\SX,\sB)$ for which both equalities are attained.
\end{theorem}

We discuss a bit these theorems and their relation with earlier estimates in the literature.
Theorem \ref{th2} is a variant of a result of Albiac and Ansorena \cite{AA1}, which for $\sB$ quasi-greedy and democratic showed that  \[
\tL\leq g^c \nu,\quad \mbox{ where }\quad g^c=\sup_{N\geq1} g^c_N\;\mand\;\nu=\sup_{N\geq1} \nu_N;\]
see \cite[Proposition 2.1.ii]{AA1}. In the unconditional case, they announced as well the bound $\bL\leq k^c\nu$ with $k^c=\sup k^c_N$ (see \cite[Remark 2.6]{AA1}), which itself improves the earlier bound $\bL\leq (k^c)^2\nu$ by Dilworth et al \cite[Theorem 2]{DKOSZ}.
Our (modest) contribution here is the explicit dependence on $N$ of the involved constants, together with a slightly shorter and more direct proof.
As discussed in \cite{AA1}, the main interest of these estimates occurs when
$\sB$ is an unconditional basis with $k_N^c\equiv 1$. Actually,  \eqref{lower}, \eqref{LN2} and the trivial estimate \[\tL\leq \bL\leq k^c_N\,\tL\] (see \cite[(1.7)]{GHO}), give

\begin{corollary}\label{CorTh2}
If for some $N$ we have $k_N^c= 1$, then
\[
\bL=\tL=\nu_N\,.
\]
\end{corollary}

In particular, the optimality asserted in the last sentence of Theorem \ref{th2}
is attained for any 1-suppression unconditional basis $\sB$.
We discuss other examples in $\S\ref{examples}$ below.

Theorem \ref{th2}, however, has some drawbacks, the first one being that in practice
$\nu_N$ may be much harder to compute explicitly than the standard
democracy constants $\mu_N$ and $\tmu_N$. A second drawback comes from the multiplicative bound $k_{2N}^c\nu_N$,
which may be far from optimal when both $k_N^c$ and $\nu_N$ grow to $\infty$.
This already occurs with simple examples of quasi-greedy bases.

\

Theorem \ref{th3} intends to cover some of these drawbacks, with an estimate which is asymptotically
optimal at least for quasi-greedy bases. In fact, if we set
\Be
\bq:=\sup_N\hg_N\,=\,\min\big\{ \,\sup_{G\in\sG}\|G\|,\;\sup_{G\in\sG}\|I-G\|\,\big\}
\label{bq}\Ee then we can show

\begin{corollary}\label{Co4}
If $\sB$ is a quasi-greedy bases and $\SK=\SR$, then
\Be
\max\{k^c_N,\mu_N\}\,\leq \,\bL\,\leq\, k^c_{2N} \,+\,8\,\bq^2\,\mu_N\label{kg2}\Ee
and
\Be
\max\{g^c_N,\mu_N\}\,\leq \,\tL\,\leq\, g_N^c \,+\,8\,\bq^2\,\mu_N\label{gg2}.\Ee
If $\SK=\SC$, the same holds with the last summand multiplied by $4$.
\end{corollary}

The fact that  $\bL\approx k_N + \mu_N$ for quasi-greedy bases is already known \cite{GHO}.
Our contribution here is an improvement of the implicit constants in the second summand, compared to $O(\bq^4)$ in \cite{GHO}, and $8\bq^3$ in \cite{DKO}. Similarly, for $\tL$ the earlier estimates in \cite[Theorem 2]{TYY2} only gave $8\bq^4$ for the involved constants in the second summand.

Another application of Theorem \ref{th3} is to bases $\sB$ which are superdemocratic but not necessarily quasi-greedy (see e.g. \cite[Example 4.8]{DKK}). In this case we have asymptotically optimal bounds $\bL\approx k_N$ and $\tL\approx g_N$; see Example \ref{Ex5} below.

\

Finally, we should say that the estimates in \eqref{LN3}, being multiplicative,
suffer from a similar drawback as \eqref{LN2}, namely they may be far from efficient when both $\tmu_N$ and $g_N$ grow fast to infinity. For such cases one always has the following trivial upper bounds

\begin{theorem}\label{th1}
If $\ttK=\sup_{m,n}\|\be_m\|\|\ben\|$, then for all $N\geq 1$ we have
\Be
\bL \leq \,1\,+\,3\ttK \;N\;\mand\;\nu_N\leq \tL \leq \,1\,+\,2\ttK \;N.
\label{LN1}\Ee
Moreover, there exists an example of $(\SX,\sB)$ for which all the equalities hold.
\end{theorem}

The optimality  for $\bL$ in Theorem \ref{th1} was first proved by Oswald \cite{Osw}.
We give a different and simpler example in $\S\ref{examples}$ below.


\section{Some elementary lemmas}
\subsection{Truncation operators}

For each $\al>0$, we define the $\al$-truncation of $z\in\SC$ by
\[
T_\al(z) = \al \,\sig(z) \;\;{\rm if}\; |z|\geq\al,\mand T_\al(z) = z \;\;{\rm if}\; |z|\leq\al.\]
We extend $T_\al$ to an operator in $\SX$ by
\Be
T_\al(x)=\sum_n\,T_\al(\ben(x))\be_n=\sum_{n\in\La_\al}\al\tfrac{\ben(x)}{|\ben(x)|}\be_n\,+\,
\sum_{n\not\in\La_\al}\ben(x)\be_n,\label{Tal}\Ee
where $\La_\al=\{n\mid|\ben(x)|>\al\}$. Since $\La_\al$ is a finite set, the last summand can be expressed
 as $(I-P_{\La_\al})x$, so the operator is well-defined for all $x\in\SX$.

\begin{lemma}\label{L2}
If $x\in\SX$ and $\bfe=\{\sig \ben(x)\}$, then \Be
\min_\La|\ben(x)|\,\big\|\bone_{\bfe\La}\big\|\,\leq\,\tilde g_N\,\|x\|,\quad \forall\;\La\in\sG(x,N).\label{abel}\Ee
\end{lemma}
\Proof
Set $\al=\min_\La|\ben(x)|$. Notice first that\Be 
T_\al x\,=\,\int_0^1\Big[\sum_n\,\chi_{[0,\frac\al{|\ben(x)|}]}(s)\,\ben(x)\be_n\Big]\,ds\,=\,
\int_0^1(I-P_{\La_{\al,s}})x\,ds,
\label{auxL1}\Ee
where we have set $\La_{\al,s}=\{n\mid|\ben(x)|>\frac\al s\}$ for each $s\in(0,1]$.

 Hence
\[
\al\bone_{\bfe\La}\,=\,
T_\al x- P_{\La^c} x\,= \int_0^1(P_\La x-P_{\La_{\al,s}}x)\,ds.\]
Note that $\La_{\al,s}\in\sG(x,k_s)$ with $k_s=|\La_{\al,s}|$ and $\La_{\al,s}\subseteq \Lambda_\alpha \subset \Lambda$. Hence
$$\|P_\La x-P_{\La_{\al,s}}x\|\le \tilde g_N \|x\|, \quad 0<s\le 1.$$
The result now follows.
\ProofEnd
\BR
The inequality \Be
\al\,\big\|\bone_{\bfe\La}\big\|\,\leq\,2\min\{g_N,g^c_N\}\,\|x\|.\label{abel1}\Ee was also  proved by an elementary Abel summation argument;
see \cite[Lemma 2.2]{DKKT}.
\ER

The next lemma is a slight improvement over \cite[Proposition 3.1]{DKK}.

\begin{lemma}
\label{L3}
For all $\al>0$, $|A|<\infty$ and $x\in\SX$ we have
\Be
\|T_\al x\|\,\leq\,g^c_{|\La_\al|}\,\|x\|,\quad \|(I-T_\al) x\|\,\leq\,g_{|\La_\al|}\,\|x\|, \label{auxL2}
\Ee
and
\Be
\|T_\al(I-P_A)x\|\,\leq\,k^c_{|A\cup\La_\al|}\,\|x\|,
\label{TalP}\Ee
where $\La_\al=\{n\mid|\ben(x)|>\al\}$.\end{lemma}
\Proof
The result follows Minkowsky's inequality and the formulae \eqref{auxL1},\[
(I-T_\al) x\,=\,
\int_0^1 \,P_{\La_{\al,s}}x\,ds.
\] and
\[
T_\al(I-P_A) x\,=\,\int_0^1(I-P_{\La_{\al,s}})(I-P_A)x\,ds,\,=\,
\int_0^1(I-P_{A\cup\La_{\al,s}})x\,ds.
\]
\ProofEnd

\BR\label{talp}
Of course, together with \eqref{TalP} one has the trivial estimate\Be
\|T_\al(I-P_A)x\|\,\leq\,g^c_{|\La_\al|}\,k^c_{|A|}\,\|x\|.
\label{TalP2}\Ee
Being multiplicative, \eqref{TalP2} is typically worse than \eqref{TalP}
(if say both $k^c_N$ and $g^c_N$ grow fast as $N\to\infty$). However in some cases it may better (e.g. when $g^c_{|\La_\al|}=1$).
\ER

\subsection{Convex extensions}

We shall use an elementary convexity lemma. As usual, the convex envelop of a set $S$ is defined by
$\co S=\{\sum_{j=1}^n\la_jx_j\mid$ {\footnotesize$x_j\in S,\;0\leq \la_j\leq1,\;\sum_{j=1}^n\la_j=1,\;n\in\SN$}$\}$.

\begin{lemma}\label{L4}
For every finite $A\subset\SN$, we have
\[
\co\Big\{\bone_{\bfe A}\mid \bfe\in\Upsilon\Big\}\,=\,\Big\{\sum_{n\in A} z_n\be_n\mid |z_n|\leq 1\Big\}.
\]\end{lemma}
\Proof We sketch the proof in the complex case, where it may be less obvious. The inclusion ``$\subseteq$''
is clear, since each $\bone_{\bfe A}$ belongs to the set $R$ on the right hand side, and $R$ is a convex set.
To show ``$\supseteq$'' one proceeds by induction in $N=|A|$. It is clear for $N=1$, so we show the case $N$ from the case $N-1$. We may assume that $A=\{\be_1,\ldots,\be_N\}$. Pick any $z=\sum_{n=1}^Nz_n\be_n\in R$, that is $|z_n|\leq1$. Write $z_N=re^{i\theta}$, and by the induction hypothesis\[
z'=\sum_{n=1}^{N-1}z_n\be_n=\sum_{\bfe}\la_\bfe \,(\e_1\be_1+\ldots+\e_{N-1}\be_{N-1}),\]
for suitable numbers $0\leq \la_\e\leq 1$ such that $\sum_\bfe \la_\e=1$. Then we have\Beas
z & = &\tfrac{1+r}2\,\big[z'+e^{i\theta}\be_N\big]\,+\,\tfrac{1-r}2\,\big[z'-e^{i\theta}\be_N\big]\\ & =&
\sum_{\bfe, \pm}\tfrac{1\pm r}2\la_\e \,(\e_1\be_1+\ldots+\e_{N-1}\be_{N-1}\pm e^{i\theta}\be_N),
\Eeas
which belongs to the set on the left hand side.
\ProofEnd

The next lemma is a straightforward extension of the inequality defining $\nu_N$.

\begin{lemma}
\label{L5}
Let $x\in\SX$ and $\al\geq\max|\ben(x)|$. Then
\[
\big\|x+z\big\|\,\leq \,\nu_N\,\big\|x\,+\,\al\bone_{\bfeta B}\big\|,\quad \forall\;\bfeta\in\Upsilon
\]
 and for all $B$ and $z$ such that $|\supp z|\leq |B|\leq N$, $B\cupdot x\cupdot z$ and $|z|_\infty\leq\al$.
\end{lemma}
\Proof We may assume that $\al=1$. By definition of $\nu_N$, the result is true when $z=\bone_{\bfe A}$,
for any $\bfe\in\Upsilon$ and any set $A$ with $|A|=|B|$ and $A\cupdot B\cupdot x$. By convexity of the norm, it continues to be true for any $z\in\co\big\{\bone_{\bfe A}\mid \bfe\in\Upsilon\big\}$. Then the general case follows
from Lemma \ref{L4}.\ProofEnd

In a similar fashion one shows
\begin{lemma}
\label{L6}
Let  $z\in\SX$ and $B\subset \SN$ such that $|\supp z|\leq|B|\leq N$. Then
\[
\big\|z\big\|\,\leq \,\tmu_N\,\max|\ben(z)|\,\big\|\bone_{\bfeta B}\big\|,\quad \forall\;\bfeta\in\Upsilon.
\]
\end{lemma}

\section{Proof of the theorems}

The general outline for proving estimates of $\bL$ and $\tL$ goes back to the work of Konyagin and Temlyakov
\cite{KT}, with the improvements coming from refinements in certain steps. In Theorem \ref{th2} we use the ideas developed by Albiac and Ansorena \cite{AA1}, slightly simplified according to our previous lemmas.

\subsection{Proof of Theorem \ref{th2}}
Let $x\in\SX$ and $\Ga\in\sG(x,N)$, and call $\al=\min_\Ga|\ben(x)|$.
Pick any $z\in\Sigma_N$ and $A\supset \supp z$ with $|A|=|\Ga|=N$.
Then we can write
\Be
x-P_\Ga x \,=\, (I-P_{A\cup\Ga})x \,+\, P_{A\setminus\Ga}x\,=:\,X+Z.
\label{first}\Ee
Since $|X|_\infty, |Z|_\infty\leq\al$ and $|\supp Z|\leq|A\setminus\Ga|=|\Ga\setminus A|$,
we can apply Lemma \ref{L5} with $\bfeta=\{\sig\ben(x)\}$ to obtain
\Bea
\|x-P_\Ga x\|& \leq & \nu_N\,\big\|\al\bone_{\bfeta(\Ga\setminus A)}\,+\,P_{(A\cup\Ga)^c}x\big\|\nonumber\\
& = & \nu_N\,\big\|T_\al\big[(I-P_A)x\big]\big\|\,=\,\nu_N\,\big\|T_\al\big[(I-P_A)(x-z)\big]\big\|\nonumber\\
& \leq & \nu_N\,k^c_{|A\cup\Ga|}\,\|x-z\|\,\leq\, \nu_N\,k^c_{2N}\,\|x-z\|,\label{auxtalp}\Eea
using Lemma \ref{L3} in the second to last inequality. Thus, taking the infimum over all $z\in\Sigma_N$
we conclude that \[\bL\leq \nu_N\,k^c_{2N}.\]

\

The estimate for $\tL$ is similar: for any set $A$ with $|A|=|\Ga|=N$ we have
\[
\|x-P_\Ga x\| \leq  \nu_N\,\big\|T_\al\big[(I-P_A)x\big]\big\|\,\leq\,\nu_N\,g^c_N\,\|x-P_Ax\|,\]
using again Lemma \ref{L3} (and $|\La_\al|\leq|\Ga|=N$). By a standard perturbation argument as in \cite[Lemma 2.2]{AA1}, this inequality continues
to hold for all $|A|\leq N$. This implies that $\tL\leq \nu_Ng^c_N$, and establishes the theorem.
\ProofEnd
\BR
Notice that we could use in \eqref{auxtalp} the estimate in Remark \ref{talp},
leading to the slightly smaller bound\[
\bL\,\leq\,\min\{k^c_{2N}, k^c_Ng^c_N\}\,\nu_N.\]
For instance, if $k^c_N=g^c_N=1$ for some $N$, this implies $\bL=\nu_N$ (as asserted in Corollary \ref{CorTh2}). In particular, one always has $\mathbf{L}_1=\nu_1$ (at least for normalized systems $\|\be_n\|=\|\ben\|=1$).
\ER

\subsection{Proof of Theorem \ref{th3}}
With the same notation as in {\eqref{first}, it is clear that\Be
\|(I-P_{A\cup\Ga})x\|\,= \,\|(I-P_{A\cup\Ga})(x-z)\|\,\leq \,k^c_{2N}\,\|x-z\|.
\label{S1}\Ee
So we only need to estimate the term $\|P_{A\setminus\Ga}x\|$.
We pick any set $\tG\in\sG(x-z,|A\setminus\Ga|)$, and use the elementary observation
\Be
\max_{A\setminus\Ga}|\ben(x)|\,\leq\,\min_{\tGa}|\ben(x-z)|;
\label{key1}\Ee
see e.g. \cite[p. 453]{GHO}. Then,
 Lemma \ref{L6} with $\bfeta=\{\sig\ben(x-z)\}$, followed by \eqref{key1} and Lemma \ref{L2} give \Bea
\big\|P_{A\setminus\Ga}x\big\| & \leq & \tmu_N\,\max_{A\setminus\Ga}|\ben(x)|\,\|\bone_{\bfeta\tGa}\|\nonumber\\
& \leq & \tmu_N\,\min_{\tGa}|\ben(x-z)|\,\|\bone_{\bfeta\tGa}\|\nonumber\\
& \leq & \tmu_N\,\tilde g_N\;\|x-z\|.\label{S2}
\Eea
So, adding up \eqref{S1} and $\eqref{S2}$ and taking the infimum over all $z\in\Sigma_N$ one obtains
\[
\|x-Gx\|\,\leq\,\big(k^c_{2N}\,+\,\tmu_N\,\tilde g_N\Big)\;\sigma_N(x),\]
as asserted in \eqref{LN3}.

\

The estimate for $\tL$ is again similar: given a set $A$ with $|A|=|\Ga|=N$,
we can replace \eqref{S1} by
\Be
\|(I-P_{A\cup\Ga})x\|\,= \,\|(I-P_{\Ga\setminus A})(I-P_A)x\|\,\leq \,g^c_{N}\,\|x-P_Ax\|,
\label{S1b}\Ee
since $\Ga\setminus A\in\sG\big(x-P_Ax\big)$. The second estimate in \eqref{S2}
is valid in this case setting $z=P_Ax$ and $\tGa=\Ga\setminus A$.
Thus we conclude
\[\|x-G_Nx\|\,\leq\,\big(g^c_{N}\,+\,\tmu_N\,\tilde g_N\,\big)\;\inf_{|A|=N}\|x-P_Ax\|,\]
and as before, this last quantity coincides with $\widetilde\sigma_N(x)$ by \cite[Lemma 2.2]{AA1}. The optimality of the constants is a consequence of Example \ref{l1c0}, that we discuss below.
\ProofEnd
\BR
In \eqref{S1} one could replace $k^c_{2N}$ by $g^c_N\,k^c_N$, arguing as in \eqref{S1b}. Typically, the latter will be a worse constant,
except in some special cases, such as if $k_N^c=1$ for some $N$, in which case $\bL=\tL\leq 1+\tmu_N$ (regardless of what $k^c_{2N}$ could be).
\ER

\subsection{Proof of Theorem \ref{th1}}

The first estimate in \eqref{LN1} is implicit in the first papers in the topic (see e.g., \cite{Tem98trig, Tem98haar} or \cite[(1.8)]{Osw}). We sketch below the elementary proof, as it also gives the second estimate.
With the notation in  \eqref{first}, notice that
\Bea
\big\|P_{A\setminus\Ga}x\big\| & \leq & \sum_{m\in A\setminus\Ga}|\be^*_m(x)|\|\be_m\|\; \leq \; \sup_m\|\be_m\|\,\sum_{n\in \Ga\setminus A}|\ben(x)|\nonumber\\
& \leq & \sup_{m,n}\|\be_m\|\|\ben\|\;N\;\|x-z\|,\label{S2b}
\Eea
since $\ben(x)=\ben(x-z)$ when $n\not\in A$. Thus, using either \eqref{S1} or \eqref{S1b}
we see that\Be
\bL\,\leq\, k_{2N}^c+\ttK \,N\;\mand \;\tL\leq g_N^c+\ttK\,N.
\label{LkN}\Ee
Now \eqref{LN1} follows from \eqref{LkN} and the trivial upper bound
\Be
k_N \,\leq\,\ttK_* \,N\;\Longrightarrow\;g^c_N\leq k^c_N \,\leq\,1\,+\,\ttK_* \,N,
\label{KN}\Ee
with  $\ttK_*=\sup_{n\geq1}\|\be_n\|\|\ben\|\leq \ttK$.
The optimality of the constants is a consequence of Example \ref{Ex1}, that we discuss below.
\ProofEnd

\subsection{Proof of Corollary \ref{Co4}}

We need an additional inequality to pass from $\tmu_N$ to $\mu_N$.
Consider the new constant
\Be
\ga_N=\sup\,\Big\{\frac{\|\bone_{\bfe B}\|}{\|\bone_{\bfe A}\|}\mid B\subset A,\;|A|\leq N,\; \bfe\in\Upsilon\Big\},
\label{gaN}
\Ee
and observe that $\ga_N\leq\hg_N$. We also have the following
\begin{lemma}
\label{LUCC}
Let $\kappa=1$ or $2$, if $\SX$ is real or complex, respectively. Then,
\Be
\big\|\bone_{\bfe B}\big\|\,\leq\,
2\kappa\,\ga_N\,\big\|\bone_{\bfeta A}\big\|,\quad \forall\;B\subset A,\;|A|\leq N,\;\bfe,\bfeta\in\Upsilon.
\label{UCC}\Ee
\end{lemma}
\Proof
Observe that changing the basis $\{\be_n\}$ to $\{\eta_n\be_n\}$ does not modify the value of $\ga_N$.
So we may assume in \eqref{UCC} that $\bfeta\equiv 1$.
We use the convexity argument in \cite[Lemma 6.4]{DKO}.
First notice that \eqref{gaN} actually implies
\Be
\|x\|\leq\ga_N\|\bone_A\|,\quad\forall\;x\in S=\Big\{\sum_{A'\subset A}\theta_{A'}\bone_{A'}\mid \sum_{A'\subset A}|\theta_{A'}|\leq 1\Big\}.
\label{Sset}\Ee
In the real case, splitting $B=B_+\cupdot B_-$, with $B_\pm=\{n\in B\mid \e_n=\pm1\}$, it is clear that
$\bone_{\bfe B}=\bone_{B_+}-\bone_{B_-}\in 2S$. In the complex case, a slightly longer argument  as in \cite[Lemma 6.4]{DKO} gives that $\bone_{\bfe B}\in 4S$.
So, in both cases we obtain \eqref{UCC}.
\ProofEnd

\begin{lemma}
\label{LtmuN}
Let $\kappa$ be as in Lemma \ref{LUCC}. Then,
\Be
\tmu_N\,\leq\, 4\,\kappa^2\,\ga_N\,\mu_N,\quad \forall\;N=1,2,\ldots
\label{tmuN}\Ee
\end{lemma}
\Proof
Take $A,B\subset\SN$  with $|A|=|B|\leq N$ and $\bfe,\bfeta\in\Upsilon$. We must show that
\Be
\|\bone_{\bfe A}\|\,\leq\, 4\,\kappa^2\,\ga_N\,\mu_N\,\|\bone_{\bfeta B}\|.\label{tmuN2}
\Ee
In the real case, split $A=A_1\cupdot A_2$ with $A_j=\{n\in A\mid \e_n=(-1)^j\}$,
and pick any partition $B=B_1\cupdot B_2$ such that $|B_j|=|A_j|$, $j=1,2$. Then
\[
\|\bone_{\bfe A}\|\,\leq\,\|\bone_{A_1}\|+\|\bone_{A_2}\|\,\leq\,
\mu_N\,\big[\|\bone_{B_1}\|+\|\bone_{B_2}\|\big]\,
\leq\,4\,\ga_N\,\mu_N\,\|\bone_{\bfeta B}\|,\]
using Lemma \ref{LUCC} in the last step. In the complex case, arguing as in \eqref{Sset} from the previous lemma, we have
 $\bone_{\bfe A}\in 4S$.
Now given  $x=\sum_{A'\subset A}\theta_{A'}\bone_{A'}\in S$,
we pick for each $A'$ a subset $B'\subset B$ such that $|A'|=|B'|$. Again, we have
\[
\|x\|\,\leq\,\sum_{A'\subset A}|\theta_{A'}|\|\bone_{A'}\|\,\leq\,\mu_N\,
\sum_{A'\subset A}|\theta_{A'}|\|\bone_{B'}\|\,\leq\,\mu_N\,
\,2\kappa\,\ga_N\,\|\bone_{\bfeta B}\|,\]
using Lemma \ref{LUCC} at the last step.
This easily gives \eqref{tmuN2}.
\ProofEnd

\Proofof{Corollary \ref{Co4}} By Theorem \ref{th3} and Lemma \ref{LtmuN}, the last summand in \eqref{LN3}
can now be controlled by
\[
\hg_N\,\min\{2,\vg_N\}\,\tmu_N\,\leq\, 2\hg_N\,4\kappa^2\,\ga_N\,\mu_N\,\leq\, 8\kappa^2\,\hg_N^2\,\mu_N.
\]
This clearly implies \eqref{kg2} and \eqref{gg2}.\ProofEnd

\BR\label{Rnew1}
Observe that we actually have the more general bounds
\Be
\bL\,\leq\,k^c_{2N}+8\kappa^2\,\ga_N\,\hg_N\,\mu_N\,,\mand \tL\,\leq\,g^c_{N}+8\kappa^2\,\ga_N\,\hg_N\,\mu_N\,.
\label{new1}\Ee
We show in Example \ref{Ex5} below that this bound is asymptotically optimal for some non quasi-greedy bases.
\ER


\section{Lower bounds: proof of Proposition \ref{P1}}

The lower bounds in \eqref{lower} are quite elementary,
and most of them have appeared before in the literature.
We sketch the proof of those we did not find explicitely in this generality.

\subsection{$\bL\geq k^c_N$} This can be found in \cite[Proposition 3.3]{GHO}.

\subsection{$\tL\geq\mu_N$} For any $|A|=|B|\leq N$, let
\[
x=\bone_{A\setminus B} + \bone_{B\setminus A} + \bone_{A\cap B}+ \bone_C,\]
where $C$ is any set such that $A\cupdot B\cupdot C$ and $|A\setminus B|+|C|=N$.
Then we can select $G_N\in\sG_N$ such that $G_Nx=\bone_{A\setminus B} +\bone_C$ and obtain
\[
\|1_B\|=\|x-G_Nx\|\leq\tL \tsigma_N(x)\leq \tL \|x-P_{C\cup B\setminus A}x\| =
\tL\,\|\bone_A\|.\]
This clearly implies $\tL\geq \mu_N$.
\BR
A similar construction can be used to show that
\[
\tL\geq \tmu^d_N = \sup\Big\{\frac{\|\bone_{\bfe A}\|}{\|\bone_{\bfeta B}\|}\mid |A|=|B|\leq N, \;A\cap B=\emptyset,\;\bfe,\bfeta\in\Upsilon\Big\}.
\]
We do not know whether one may actually have $\tL$ or even $\bL\geq \tmu_N$.\ER

\subsection{$\tL\geq \frac1{2\kappa}\tmu_N$}
Given $|A|=|B|\leq N$ and $\bfe,\bfeta\in\Upsilon$, we must show that
\[
\|\bone_{\bfeta B}\|\,\leq \,2\kappa\, \tL\,\|\bone_{\bfe A}\|.
\]
It is enough to prove it for $\bfe\equiv1$ (otherwise, apply the result to $\sB=\{\e_n\be_n\}$).
Recall from \eqref{Sset} (and \cite[Lemma 6.4]{DKO}) that $\bone_{\bfeta B}\in 2\kappa S$, where
\[S=\Big\{\sum_{B'\subset B}\theta_{B'}\bone_{B'}\mid \sum_{B'\subset B}|\theta_{B'}|\leq 1\Big\},
\]
so it suffices to show that
\[
\|\bone_{B'}\|\,\leq \, \tL\,\|\bone_{A}\|,\quad\forall\;B'\subset B.
\]
Pick any $C\subset (A\cup B)^c$ with $|A\setminus B'|+|C|=N$ and set
\[
x=\bone_{B'\setminus A}+\bone_{B'\cap A}+\bone_{A\setminus B'}+\bone_C.\]
Then can take $G_N\in\sG_N$ such that $G_Nx=\bone_{A\setminus B'}+\bone_C$, and hence
\[
\|\bone_{B'}\|=\|x-G_Nx\|\leq\tL\tsigma_N(x)\leq \tL\|x-P_{C\cup(B'\setminus A)}x\|=\tL\|\bone_A\|,
\]
where we have used $|B'\setminus A|\leq|B\setminus A|=|A\setminus B|\leq |A\setminus B'|=N-|C|$.

\subsection{$\tL\geq\nu_N$}
Let $|A|=|B|\leq N$, $\bfe,\bfeta\in\Upsilon$, and $x\in\SX$ such that $A\cupdot B\cupdot x$ and $|x|_\infty\leq 1$.
We must show that
\Be
\|\bone_{\bfe A} + x\| \,  \leq \, \tL \,\|\bone_{\bfeta B} + x\|,
\label{auxx2}\Ee
For every $j\geq1$ we can find a set $C_j$ with $|C_j|=N-|A|$, disjoint with $A\cup B$,  and  such that $\max_{n\in C_j}|\ben(x)|\leq 1/j$. We set
\[
y_j=\bone_{\bfe A} + \bone_{\bfeta B} + (I-P_{C_j})x+ \bone_{C_j},\]
and select $G_N\in\sG_N$ such that $G_N(y_j)=\bone_{\bfeta B} + \bone_{C_j}$. Then
\Beas
\|\bone_{\bfe A} + (I-P_{C_j})x\| & = & \|y_j - G_N(y_j)\| \,  \leq \, \tL \,\tsigma_N(y_j) \\
& \leq & \tL\,
\|(I-P_{A\cup C_j})y_j\| =\tL\|\bone_{\bfeta B} + (I-P_{C_j})x\|.\Eeas
Since  $\lim_{j\to\infty}P_{C_j}x=0$ we obtain \eqref{auxx2}.

\subsection{$\tL\geq g^c_N$} We must show that for every $x\in\SX$ and every $\Ga\in\sG(x,k)$ with $k\leq N$, we have
\Be
\|x-P_\Ga x\|\leq\tL\|x\|.\label{auxx1}
\Ee
Let $\al=\min_{n\in\Ga}|\ben(x)|$. Notice that for every $j\geq1$ we can find a set $C_j\subset \Ga^c$,
with $|C_j|=N-k$, and $\max_{n\in C_j}|\ben(x)|\leq \al/j$. Let
\[
y_j=x-P_{C_j}x+\al\bone_{C_j},\]
so that $\Ga\cupdot C_j\in\sG(y_j,N)$. Thus\[
\|y_j-P_{\Ga\cup C_j}y_j\|\leq \tL\,\tsigma_N(y_j)\,\leq\,\tL\,\|y_j-P_{C_j}y_j\|,\]
which is the same as \[
\|x-P_\Ga x - P_{C_j}x\|\,\leq\,\tL\,\|x-P_{C_j}x\|.
\]
Since $\lim_{j\to\infty}P_{C_j}x=0$ (in $\SX$) we obtain \eqref{auxx1}.
\ProofEnd

\section{Examples}\label{examples}

\subsection{The summing basis}\label{Ex1}
Let $\SX$ be the (real) Banach space of all sequences $\ba=(a_n)_{n\in\SN}$ with
\Be
\|\ba\|:=\sup_{M\geq1}\Big|\sum_{n=1}^M a_n\Big|<\infty.\label{Ex1norm}\Ee
The standard canonical basis $\{\be_n,\ben\}$ satisfies $\|\be_m\|\equiv1$, $\|\be^*_1\|=1$ and $\|\ben\|=2$ if $n\geq2$ (so $\ttK=2$, with the notation in Theorem \ref{th1}). The terminology comes from the fact that $\SX$ is isometrically isomorphic\footnote{Via the map $\ba\in\SX\mapsto T\ba=(\sum_{i=1}^n a_i)_{n\in\SN}\in\ell^\infty$, since $T\be_n=\bs_n$.} to the span of the ''summing system'' $\{\bs_n:=\sum_{k\geq n}\be_k\}_{n=1}^\infty$ in $\ell^\infty$; see \cite[p. 20]{LZ}.

\begin{proposition}
For this example we have\Bi\item $\mu_N=1$ and $\tmu_N=N$
\item $g_N=k_N=2N$ and $g^c_N=k^c_N=1+ 2N$
\item $\nu_N=\tL=1+ 4N$  and $\bL=1+6N$.\Ei
So, equalities hold everywhere in Theorem \ref{th1}.
\end{proposition}
\Proof
It is clear that $\|\bone_A\|=|A|$, so the basis is democratic and $\mu_N\equiv1$. On the other hand, we trivially have \[
1\leq \|\bone_{\bfe A}\|\leq N, \quad \forall\;|A|=N, \;\bfe\in\Upsilon.
\] The upper bound is attained if $\bfe\equiv1$, and the lower bound is attained in the explicit example
$\|\sum_{n=1}^N(-1)^n\be_n\|=1$. We conclude that $\tmu_N=N$.

\

We know from \eqref{KN} that $g_N\leq k_N\leq 2N$. To see the equality, pick the vector
$\ba 
=(-1,2,-2,\ldots,2,-2,0,\ldots)$, which has $\|\ba\|=1$. Then
$\Ga=\{n\!\!\mid\!\! a_n=2\}\in\sG(\ba,N)$ and
\[
g_N\geq \|P_\Ga\ba\|=\|(0,2,0,\ldots,2,0,0\ldots)\|=2N.\]
Similarly, $g_N^c\leq k_N^c\leq 1+2N$ by \eqref{KN}, and
setting $\Ga'=\{n\!\!\mid\!\! a_n=-2\}\in\sG(\ba,N)$ we conclude
\[
g^c_N\geq \|(I-P_{\Ga'})\ba\|=\|(1,2,0,\ldots,2,0,0\ldots)\|=1+2N.\]
Next we have $\nu_N\leq \tL\leq 1+4N$, by Proposition \ref{P1} and Theorem \ref{th1}. For the lower bound we pick
\[
x=\big(\overbrace{\tfrac12,0,\tfrac12}\;;\;\ldots\;;\;\overbrace{\tfrac12,0,\tfrac12}\;; \;\tfrac12,0,0,\ldots\big)\mand \bone_B= \big(\overbrace{0,1,0}\;;\;\ldots\;;\;\overbrace{0,1,0}\;; \;0,\ldots\big)
\]
so that $\|x-\bone_B\|=1/2$, while $\|x+\bone_A\|=\frac12+2N$ for any $|A|=N$.
So,
\[
\nu_N\geq \frac{\|x+\bone_A\|}{\|x-\bone_B\|}\,=\,1+4N.\]
Finally, $\bL\leq 1+6N$ by Theorem \ref{th1}. To show equality, let
\[
x=\big(\overbrace{\tfrac12,1,\tfrac12}\;;\;\ldots\;;\;\overbrace{\tfrac12,1,\tfrac12}\;; \;\tfrac12;\;
\overbrace{-1,1},\;\ldots\;,\overbrace{-1,1}\,,\;0,0,\ldots\big),
\]
and pick $\Ga=\{n\!\!\mid\!\! x_n=-1\}\in\sG(x,N)$. Then
\[
\|x-P_\Ga x\|\,=
\,3N+\tfrac12,
\]
while
\Beas
\sigma_N(x) & \leq & \big\|x-2\big(\overbrace{0,1,0}\;;\;\ldots\;;\;\overbrace{0,1,0}\;;\;0,0,\ldots\big)\big\|
\,=\,\frac12.\\
\Eeas
Thus, $\bL\geq {\|x-P_\Ga x\|}/{\sigma_N(x)}\,\geq\,6N+1$.
\ProofEnd

\BR
In this example one can also show that $\ga_N=\lceil N/2\rceil$ for the constant defined in \eqref{gaN}. In particular, the bound in \eqref{UCC} (with $\kappa=1$) cannot be improved.
\ER

\subsection{Canonical basis in $\ell^1\oplus c_0$}\label{l1c0} That is, we consider pairs of sequences $(x,y)\in \ell^1\times c_0$, endowed with the norm $\|(x,y)\|=\|x\|_1+\|y\|_\infty$. Write the canonical basis as $\sB=\{(\be_m,0),(0,\mathbf{f}_n)\}_{m,n=1}^\infty$.
\begin{proposition}
The canonical basis in $\ell^1\oplus c_0$ satisfies\Bi\item $\mu_N=\tmu_N=N$
\item $g_N=k_N=g^c_N=k^c_N=1$
\item $\nu_N=\tL=\bL=1+\tmu_N=1+N$.\Ei
So, equalities hold everywhere in Theorems \ref{th2} and \ref{th3}.
\end{proposition}
\Proof
The second point is clear, since the canonical basis is 1-unconditional.
For the first point just notice that\[
1\leq \|\bone_A\|=\|\bone_{\bfe A}\|\leq |A|,\]
with the lower bound attained when $\bone_A\in c_0$, and the upper bound when  $\bone_A\in\ell^1$.
Finally, in view of Theorem \ref{th3} and the previous equalities, in the last point we only need to show that $\nu_N\geq N+1$.
Let $\bone_A=\sum_{n=1}^N\be_n$, $\bone_B=\sum_{n=1}^N\mathbf{f}_n$, and $x=\mathbf{f}_{N+1}$, then
\[
\nu_N\geq\,\frac{\|\bone_A+x\|}{\|\bone_B+x\|}=N+1.\]
\ProofEnd
\subsection{Canonical basis in $\ell^1\oplus \ell^q$, $1\leq q<\infty$}
This variant of the previous example also admits explicit Lebesgue constants, but equality fails in \eqref{LN3}.
\begin{proposition}
The canonical basis in $\ell^1\oplus \ell^q$, $1\leq q<\infty$ satisfies\Bi\item $\mu_N=\tmu_N=N^{1/q'}$
\item $g_N=k_N=g^c_N=k^c_N=1$
\item $\nu_N=\tL=\bL=(N+1)^{1/q'}$.\Ei
\end{proposition}
\Proof
We only prove the last part, the other two being easy.
By Corollary \ref{CorTh2}, we only need to estimate $\nu_N$.
From below, we choose as before $\bone_A=\sum_{n=1}^N\be_n$, $\bone_B=\sum_{n=2}^{N+1}\mathbf{f}_n$, and $x=\mathbf{f}_{1}$, so that
\[
\nu_N\geq \,\frac{\|\bone_A+\mathbf{f_1}\|}
{\|\bone_B+\mathbf{f_1}\|}=\frac{N+1}{(N+1)^\frac1{q}}\,=\,(N+1)^{1/q'}.
\]
From above, let $|A|=|B|=N$ and $(x,y)$ have disjoint support with $A\cup B$. Then\[
\|(x,y)+\bone_{\bfe A}\|\leq \|x\|_1+\|y\|_q+N,
\]
while if  
$k=|\supp P_{\ell^1}(\bone_B)|$, then\[
\|(x,y)+\bone_{\bfeta B}\|\,=\, \|x\|_1 +k\,+(\|y\|_q^q+N-k)^\frac1q\,
\geq\,  \|x\|_1+(\|y\|_q^q+N)^\frac1q.
\]
So,
\[
\frac{\|(x,y)+\bone_{\bfe A}\|}{\|(x,y)+\bone_{\bfeta B}\|}\,\leq\, \frac{\|x\|_1+\|y\|_q+N}{\|x\|_1+(\|y\|_q^q+N)^\frac1q}\,\leq\, \frac{\|y\|_q+N}{(\|y\|_q^q+N)^\frac1q},\]
and the latter is easily seen to be maximized at $\|y\|_q=1$.
So $\nu_N\leq (1+N)^\frac1{q'}$, as asserted.
\ProofEnd

\BR
With similar (but slightly more tedious) computations one can show that, for $\ell^p+c_0$, $1< p<\infty$, one has
\[
\nu_N=\tL=\bL= 1+N^\frac1p,
\]
while $\tmu_N=\mu_N=1+(N-1)^\frac1p$, so again equality fails in \eqref{LN3}.
\ER

\subsection{The trigonometric system}\label{Ex4}
Consider $\sB=\{e^{i nx}\}_{n\in\SZ}$ in $L^p(\ST)$, $1\leq p\leq\infty$.
In this case, neither \eqref{LN2} nor \eqref{LN3} give good estimates, even asymptotically. By a more direct approach, Temlyakov \cite{Tem98trig} showed the following\[
c_pN^{|\frac1p-\frac12|}\leq \bL\leq 1+3N^{|\frac1p-\frac12|},
\]
for some $c_p>0$. More precisely, the following inequalities hold (if $p>1$)
\Be
c_pN^{|\frac1p-\frac12|}\leq \ga_N\leq g^c_N\leq k^c_N\leq 1+N^{|\frac1p-\frac12|},\label{gaNLp}\Ee
and
\Be c_p N^{|\frac1p-\frac12|}\leq \mu_N\leq\tmu_N=\tmu^{d}_N
\leq\nu_N\leq\tL\leq \bL\leq 1+3N^{|\frac1p-\frac12|}.\label{muNLp}\Ee
So all the involved constants have the same order of magnitude $N^{|\frac1p-\frac12|}$.
For the upper bounds in \eqref{gaNLp} and \eqref{muNLp}, see \cite[Lemma 2.1 and Theorem 2.1]{Tem98trig}.
The lower bounds are implicit in \cite[Remark 2]{Tem98trig}; for instance if $1<p\leq 2$ and $N\in2\SN$ then\Be
\mu_{N+1}\geq \frac{\|\bone_{\{1,2,\ldots,2^N\}}\|_p}{\|\bone_{\{-N/2,\ldots,N/2\}}\|_p}\,\geq\,c_p\,\frac{\sqrt N}{N^{1-\frac1p}}\,=\,c_p\,N^{\frac1p-\frac12},
\label{low1}\Ee
since the Dirichlet kernel has norm $\|D_{N/2}\|_p\approx N^{1-\frac1p}$. Likewise, by \eqref{UCC}
\Be
\ga_{N+1}\geq \tfrac14\,\frac{\|\bone_{\bfe\{-N/2,\ldots,N/2\}}\|_p}{\|\bone_{\{-N/2,\ldots,N/2\}}\|_p}\,\geq\,c'_p\,\frac{\sqrt N}{N^{1-\frac1p}}\,=\,c'_p\,N^{\frac1p-\frac12},
\label{low2}\Ee
choosing in $\bfe$ the signs of the corresponding Rudin-Shapiro polynomial. The case $p\geq2$ is similar, replacing the roles of numerator and denominator.

When $p=1$ the arguments in \cite{Tem98trig} still give
\Be
\bL\approx\tL\approx k_N\approx g_N\approx \sqrt N,
\label{LNL1}\Ee whereas \Be
\ga_N\approx\mu_N\approx\tmu_N\approx \frac{\sqrt N}{\log N}.\label{gaNL1}\Ee
In this last estimate the lower bound for each of the constants follows as in \eqref{low1} and \eqref{low2}, using $\|D_{N/2}\|_1\approx\log N$. The upper bound relies on $\|\bone_{\bfeta B}\|_1\leq\|\bone_{\bfeta B}\|_2=|B|^{\frac12}$, and on the deeper result $\inf_{\bfe, |A|=N}\|\bone_{\bfe A}\|_1\geq c\log N$, a famous problem posed by Littlewood and solved by Konyagin \cite{Kon} and McGeehee-Pigno-Smith \cite{MPS}.
Finally, we show that in this case we have
\Be
\nu_N\approx\sqrt N.
\label{nuNL1}\Ee
Since $\nu_N\leq \bL\lesssim\sqrt N$, we only need to show the lower bound.
For $N\in\SN$ we pick $B=\{-N,\ldots,N\}$ and $x$ so that\[
\bone_{\{-N,\ldots,N\}}+x=V_{N},
\]
where $V_{N}$ denotes the de la Vall\'ee-Poussin kernel (as in \cite[p. 114]{StSha}). Then
$|x|_\infty\leq 1$, $\supp x\subset\{N<|k|<2N\}$ and we have
\[
\|\bone_B+x\|_1=\|V_N\|_1\leq 3.
\] Next we pick $A=\{2^j\!\mid\! j_0\leq j\leq j_0+2N\}$
where we choose $2^{j_0}\geq4N$. Then $(I-V_{2N})(\bone_A+x)=\bone_A$, and therefore
\[
c_1\sqrt{N}\leq\|\bone_A\|_1\leq\|I-V_{2N}\|_1\,\|\bone_A+x\|_1\leq 4\,\|\bone_A+x\|_1.\]
Overall we conclude that
\[
\nu_{2N+1}\geq \frac{\|\bone_A+x\|_1}{\|\bone_B+x\|_1}\geq \tfrac{c_1}{12}\,\sqrt{N}.
\]

\subsection{A superdemocratic and not quasi-greedy basis} \label{Ex5}
Theorem \ref{th3} becomes asymptotically optimal when $\tmu_N\approx1$, as in this case $\bL\approx k_N$ and $\tL\approx g_N$. We give a non-trivial example of this situation, which is
a small variation of \cite[Example 4.8]{DKK}.
This example has the additional interesting property of being \emph{unconditional with constant coefficients}\footnote{That is, $\|\bone_{\bfe A}\|\approx\|\bone_A\|$ for all finite $A$ and all $\bfe\in\Upsilon$; see \cite[Def 3]{Wo}.} but not quasi-greedy.

\begin{proposition}\label{PEx5}
For every $1\leq q\leq\infty$, there exists $(\SX,\sB)$ such that
\Bi \item $\nu_N\approx\tmu_N\approx \ga_N\approx1$
\item $ g_N\approx k_N\approx (\log N)^{1/q'}$
\item $\bL\approx \tL\approx (\log N)^{1/q'} $
\Ei
So, in this case Theorems \ref{th2}, \ref{th3} and Remark \ref{Rnew1} are asymptotically optimal.
\end{proposition}

\Proof Let $\cD_k$ denote the set of all dyadic intervals $I\subset [0,1]$ with length $|I|=2^{-k}$,
and $\cD=\cup_{k\geq0}\cD_k$. Consider the space $\fpq$ of all (real) sequences $\ba=(a_I)_{I\in\cD}$
such that\[
\|\ba\|_{\fpq}\,=\,\Big\|\big[\sum_{I}|a_I\chii|^q\big]^\frac1q\Big\|_{L^1}<\infty,
\]
where $\chii=|I|^{-1}\chi_I$. It is well known that $\{\be_I\}_{I\in\cD}$, the canonical basis, is unconditional and democratic in $\fpq$;
see e.g. \cite{HJLY, GH}.
In particular, for some $c_q\geq1$ we have
\[
\tfrac1{c_q}|A|\leq \|\bone_{\bfe A}\|_{\fpq}\,\leq\, |A|,\quad \forall\;A\subset\cD,\;\;\bfe\in\Upsilon.\]
From the definition we also have
\[
\Big\|\sum_{k} b_k 2^{-k}\bone_{\cD_k}\Big\|_{\fpq}\,=\,\big(\sum_{k}|b_k|^q\big)^\frac1q,
\]
since $2^{-k}\sum_{\cD_k}\chii=\chi_{[0,1]}$.
For every $N\geq 1$ we shall pick a subset $\{k_1,\ldots k_N\}\subset\SN_0$, and look at the
finite dimensional space $F_N$ consisting of sequences supported in $\cup_{j=1}^N\cD_{k_j}$.
We order the canonical basis by $\cup_{j=1}^N\{\be_I\}_{I\in \cD_{k_j}}$, so we may as well write their elements
as $\ba=(a_j)_{j=1}^{d_N}$. We also consider in $F_N$ the James norm\[
\|(a_j)\|_{J_q}\,=\,\sup_{m_0=0<m_1<\ldots}\Big[\sum_{k\geq0}\big|\sum_{m_k<j\leq m_{k+1}} a_j\big|^q\Big]^\frac1q.
\]
Note that $\|\ba\|_{J_q}\leq \|\ba\|_{\ell^1}$, with equality iff all the $a_j$'s have the same sign\footnote{Note that $|a-b|<(a^q+b^q)^\frac1q$ if $a,b>0$, so consecutive elements with different signs should be in different blocks of the James norm.}. In particular,\[
\|\bone_A\|_{J_q}=|A|.\]
Now set in $F_N$ a new norm
\[
\tri{\ba}=\max\Big\{\|\ba\|_{\fpq},\;\|\ba\|_{J_q}\Big\},
\]
and observe that $1/c_q|A|\leq \tri{\bone_{\bfe A}}\leq |A|$, with $c_q$ independent of $N$ and $k_j$.
Also, the vector $x=\sum_{j=1}^N (-1)^{j+1} 2^{-k_j}\bone_{\cD_{k_j}}$ has
\[ \|x\|_{\fpq}=\|x\|_{J_q}=\tri{x}=N^\frac1q.\]
At this point we write $N=2n$ and choose our $k_j$'s as \[
k_{2j+1}=j\mand k_{2j+2}=n+j,\quad j=0,\ldots, n-1.
\]
Then if $P=\sum_{j {\rm\; odd}}2^{k_j}=2^n-1$ we have $G_Px=\sum_{j{\rm \; odd}} 2^{-k_j}\bone_{\cD_{k_j}}$, which implies
\[
\|G_Px\|_{\fpq}=n^\frac1q,\quad \|G_Px\|_{J_q}=n,\mand \tri{G_Px}=n.
\]
Therefore
\[
g_{2^n}\geq {\tri{G_Px}}/{\tri{x}}\,\geq \;n^{1-\frac1q}.
\]
We turn to estimate the unconditionality constant $k_m$ of the space $F_N$.
Given $|A|=m$,  we first claim that
\Be
 \|P_Ax\|_{\ell^1}\,\leq\,c'_q\, (\log|A|)^{1/q'}\,\|x\|_{\fpq}.
\label{Ex5aux1}\Ee
This is clear when $q=1$ (since $\mathfrak{f}_1^1=\ell^1$).
When $q=\infty$, it is a consequence e.g. of  \cite[Remark 5.6]{GH} (since $\mathfrak{f}^\infty_1$ is a $1$-space, in the terminology of \cite[(2.8)]{GH}). Thus one derives \eqref{Ex5aux1} by complex interpolation.
From here\[
\tri{P_Ax}\leq  \|P_Ax\|_{\ell^1}\,\leq\,c'_q\, (\log|A|)^{1/q'}\,\tri{x},
\]
which implies the bound $k_m\leq c_q' (\log m)^{1/q'}$.

Finally, we consider the space $\SX=\oplus_{\ell^1} F_N$ with $\sB$ the consecutive union of the natural bases in $F_N$.
Then
\[
\tfrac1{c_q}|A|\leq\tri{\bone_{\bfe A}}=\sum_N\tri{\bone_{\bfe A_N}}\leq |A|,\]
so $\sB$ is superdemocratic. We claim further that $\nu_N=O(1)$. Let $|A|=|B|=N$ and $x\in\SX$ have disjoint support with $A\cupdot B$.
Assuming first that $\tri{x}\geq 2N$, we have\[
\frac{\tri{\bone_{\bfe A}+x}}{\tri{\bone_{\bfeta B}+x}}\,\leq\,
\frac{\tri{\bone_{\bfe A}}+\tri{x}}{\tri{x}-\tri{\bone_{\bfeta B}}}\,\leq\,\frac{3/2\tri{x}}{1/2\tri{x}}=3,
\]
since $\tri{\bone_{\bfe A}},\tri{\bone_{\bfeta B}}\leq N\leq\tri{x}/2$.
Otherwise we have $\tri{x}\leq 2N$, which implies\[
\frac{\tri{\bone_{\bfe A}+x}}{\tri{\bone_{\bfeta B}+x}}\,\leq\,
\frac{\tri{\bone_{\bfe A}}+\tri{x}}{\sum_N\|\bone_{\bfeta B_N}+x_N\|_{\fpq}}\,\leq\,
\frac{3N}{\sum_N\|\bone_{\bfeta B_N}\|_{\fpq}}\,\leq\,3c_q,
\]
since $\sum_N\|\bone_{\bfeta B_N}\|_{\fpq}\geq c_q\sum_N|B_N|=N$. Thus $\nu_N\lesssim 1$ as asserted.
A similar argument shows that
\[
\ga_N\leq \frac{\tri{\bone_{\bfe A}}}{\tri{\bone_{\bfeta B}}}\,\leq\,
\frac{N}{\sum_N\|\bone_{\bfeta B_N}\|_{\fpq}}\,\leq\,c_q.
\]
Finally, observe that $k^\SX_m\leq \max_N k_m^{F_N}\leq c'_q (\log m)^{1/q'}$,
while if $N=2n$ we have \[
g_{2^n}^\SX\,\geq\, g^{F_N}_{2^n}\geq n^{1/q'}.\]
This completes the proof of Proposition \ref{PEx5}.\ProofEnd

\section{Further questions}

As shown in Example \ref{Ex4}, the multiplicative bounds in Theorems \ref{th2} and \ref{th3} are not so good when  both $g_N$ and $\tmu_N$ go to infinity.

\bline {\bf Q1:} \emph{Find bounds for $\bL$ and $\tL$ which depend
\textbf{additively} on $k_N$, $\tmu_N$ or $\nu_N$. More precisely, determine in what cases it can be true that\[
\bL \lesssim\,k_N+\nu_N\quad \mbox{ or }\quad \bL \lesssim\,k_N+\tmu_N.\]}

This is for instance the case for the trigonometric system, and the other examples in $\S\ref{examples}$.
In this respect, we can mention the results of Oswald \cite{Osw}, who obtains additive estimates of the form $\bL\approx k_N + B_N$, but with constants $B_N$ of a  more complicated nature.

\

Related to the previous one can ask

\bline {\bf Q2:} \emph{Find examples such that $k_N$ and $\nu_N$ grow independently to infinity.}

\sline Example \ref{Ex5} shows that one can have $\nu_N\approx 1$ and $\bL\approx k_N\to\infty$.
We do not know whether it is possible to have $\nu_N\approx N^\al$ and $k_N\approx N^\beta$
for arbitrary $0<\al,\beta\leq 1$.

\

The new constant $\ga_N$ in \eqref{gaN} is a natural replacement for $g_N$ in some situations.
Example \ref{Ex5} (and also \eqref{gaNL1} in Example \ref{Ex4}) show that this improvement may be strict and the ratio $g_N/\ga_N$ as large as $\log N$.

\bline {\bf Q3:} \emph{Find examples with $\ga_N\approx1$ and $g_N$ as large as possible}.

%
%
%
%
%
\bline{\bf Acknowledgements:} we wish to thank F. Albiac, J.L. Ansorena and E. Hern\'andez for many useful conversations about these topics.

\bibliographystyle{plain}

\vskip 1truemm

\end{document}